\documentclass[a4paper,10pt]{article}
\usepackage[utf8]{inputenc}
\usepackage{color}

\usepackage[all]{xy}
\usepackage{amsmath,amssymb,amsfonts,amsthm,graphicx}
\usepackage[margin=1.2in]{geometry}
\usepackage{multirow}

\newtheorem{thm}{Theorem}[section]

\newtheorem{cor}[thm]{Corollary}
\newtheorem{prop}[thm]{Proposition}

\newtheorem{lm}[thm]{Lemma}

\newtheorem*{theorem*}{Theorem}
\newtheorem*{aim*}{Aim}
\newtheorem*{initialaim*}{Initial Aim}
\newtheorem*{conj*}{Conjecture}
\newtheorem*{cor*}{Corollary}
\newtheorem*{prop*}{Proposition}
\newtheorem*{df*}{Definition}
\newtheorem*{lm*}{Lemma}
\newtheorem*{example*}{Example}
\newtheorem*{notation*}{Notation}
\newtheorem*{prob*}{Problem}

\title{On vanishing criteria that control finite group structure II}
\author{Julian Brough$^1$ and Qingjun Kong$^2$}

\begin{document}
\date{}
\maketitle

\vspace{-10mm}
\begin{center}

\small
\textit{$^1$ FB Mathematik, TU Kaiserslautern, Postfach 3049, 67653 Kaiserslautern, Germany}

\text{E-mail: brough@mathematik.uni-kl.de}

\textit{$^2$ Department of Mathematics, Tianjin Polytechnic University,
Tianjin 300387, People's Republic of China }

\text{E-mail: kqj2929@163.com}
\end{center}

\paragraph{}

\textit{Seconday: 20D06, 20D08}

\paragraph{}

\normalsize
\begin{abstract}
In a paper by the first author it was shown that for certain arithmetical results on conjugacy class sizes it is enough to only consider the vanishing conjugacy class sizes.
In this paper we further weaken the conditions to consider only vanishing elements of prime power order.
\end{abstract}

\section{Introduction}

Many results have been proven which connect the structure of a finite group $G$ to arithmetical data connected to $G$.
One type of data that has has often been considered is the set of conjugacy class sizes in a group.
Moreover, recently, instead of considering all conjugacy class sizes authors have refined this set by using the irreducible characters of a group.
In particular, the set of vanishing conjugacy class sizes has become an interesting topic. (See \cite{Brough7}, \cite{Brough3}, \cite{Brough4}, \cite{DPSVan} and also \cite{DPSVanPriGraph} for properties related to vanishing elements)
Furthermore \cite{Brough7}, \cite{Brough3}, \cite{DPSVan} and  \cite{DPSVanPriGraph} show that arithmetical data for conjugacy class sizes can be weakened to only the vanishing conjugacy class sizes.
An element $x\in G$ is called a vanishing element if there exists some irreducible character $\chi$ of $G$ such that $\chi(x)=0$; the conjugacy class $x^G$ is called a vanishing conjugacy class.

In \cite{Brough3}, the first author showed that the criterion given Cossey and Wang to determine solubility and super solubility only required the vanishing conjugacy class sizes.
Furthermore the author weakened the vanishing criterion for $p$-nilpotence given by Dolfi, Pacifici and Sanus \cite{DPSVan} to only considering the vanishing $p'$-elements (i.e. the elements whose order is not divisible by $p$).
We restate here the three main theorems given in \cite{Brough3}.

\begin{thm}\cite{Brough3}[Theorem A]
 Let $G$ be a finite group and $p$ a prime divisor of $G$ such that if $q$ is any prime divisor of $G$, then $q$ does not divide $p-1$.
 Suppose that no vanishing conjugacy class size of $G$ is divisible by $p^2$. Then $G$ is a soluble group.
\end{thm}
\begin{thm}\cite{Brough3}[Theorem B]
 Let $G$ be a finite group and suppose that every vanishing conjugacy class size of $G$ is square free.
 Then $G$ is a super soluble group.
\end{thm}

\begin{thm}\cite{Brough3}[Theorem C]
 Let $G$ be a finite group and suppose a prime $p$ does not divide the size of any vanishing conjugacy class size $|x^G|$ for $x$ a $p'$-element of prime power order in $G$.
 Then $G$ has a normal $p$-complement.
\end{thm}

The aim of this paper is to further refine which vanishing conjugacy classes are required.
In particular, it is shown that it is sufficient to only consider vanishing elements of prime power order.
In other words, we shall prove the following results:

\begin{theorem*}[Theorem 1]
 Let $G$ be a finite group and $p$ a prime divisor of $G$ such that if $q$ is any prime divisor of $G$, then $q$ does not divide $p-1$.
 Suppose that no conjugacy class size of a vanishing element of prime power order in $G$ is divisible by $p^2$. 
 Then $G$ is a soluble group.
\end{theorem*}
\begin{theorem*}[Theorem 2]
 Let $G$ be a finite group and suppose that no conjugacy class size of a vanishing element of prime power order in $G$ is square free.
 Then $G$ is a super soluble group.
\end{theorem*}

\begin{theorem*}[Theorem 3]
 Let $G$ be a finite group and suppose a prime $p$ does not divide the size of any vanishing conjugacy class size $|x^G|$ for $x$ a $p'$-element of prime power order in $G$.
 Then $G$ has a normal $p$-complement.
\end{theorem*}

Note that Theorem 1 has the following immediate corollary.

\begin{cor}
Let $G$ be a finite group and suppose that no vanishing conjugacy class size of an element of prime power order in $G$ is divisible by $4$. 
Then $G$ is a soluble group.
\end{cor}

The proof of these theorems uses very similar arguments as in \cite{Brough3}.
Therefore some of the details will be omitted here and we will instead refer the reader to the previous paper.
The key difference is that we now need to ensure our chosen elements from the non-abelian simple groups without an irreducible character of $q$-defect zero, have prime power order.
In particular, the adapted version of \cite[Lemma 2.4]{Brough3} is given by two lemmas at the end of section 2, which split the cases for the sporadic and alternating groups into two parts.

\section{Preliminaries}

Given a normal subgroup $N$ in $G$, there is a natural bijection between the set of irreducible characters of $G/N$ and the set of irreducible characters of $G$ with $N$ in their kernel.
In particular, this natural bijection implies that if $x$ is an element not in $N$ then $xN$ is vanishing in $G/N$ if and only if $x$ is vanishing in $G$.
In addition, recall that for an element $x$ in $G$, both $|x^N|$ and $|xN^{G/N}|$ divide $|x^G|$.

Let $q$ be a prime number, and $\chi$ an irreducible character of $G$; the character $\chi$ is said to have $q$-defect zero if $q$ does not divide $|G|/\chi(1)$.
A result of Brauer highlights the significance $q$-defect zero has for vanishing elements.
If $\chi$ is an irreducible character of $G$ with $q$-defect zero, then $\chi(g)=0$ for every $g\in G$ such that $q$ divides the order of $g$ \cite[Theorem 8.17]{IsaChTh}.

\begin{cor}\cite[Corollary 2]{GODefectZero}\label{ListSimpleCases1}
Let $S$ be a non-abelian simple group and assume there exists a prime $q$ such that $S$ does not have an irreducible character of $q$-defect zero.
Then $q=2$ or $3$ and $S$ is isomorphic either to one of the following sporadic simple groups $M_{12}$, $M_{22}$, $M_{24}$, $J_2$, $HS$, $Suz$, $Ru$, $Co_1$, $Co_3$, $BM$ or some alternating group $Alt(n)$ with $n\geq 7$.
\end{cor}

In the particular case that $M$ is a minimal normal subgroup, we shall use the preceding corollary together with the following lemma; this result forms a generalisation of a comment made during the proof of \cite[Theorem A]{DPSVan}.

\begin{lm}\cite[Lemma 2.2]{Brough3}\label{DefectMinLiftVan}
 Let $G$ be a group, and $N$ a normal subgroup of $G$.
 If $N$ has an irreducible character of $q$-defect zero, then every element of $N$ of order divisible by $q$ is a vanishing element in $G$.
\end{lm}

It still remains to consider those simple groups which have no character of $q$-defect zero for some prime $q$.
The next result provides a condition for an irreducible character of a minimal normal subgroup $M$ of $G$ to extend to an irreducible character of $G$.

\begin{prop}\cite[Lemma 5]{BianChDeg}\label{ExtendingChar}
 Let $G$ be a group, and $M=S_1\times \dots\times S_k$ a minimal normal subgroup of $G$, where every $S_i$ is isomorphic to a non-abelian simple group $S$.
 If $\theta\in Irr(S)$ extends to $Aut(S)$, then $\theta\times\dots\times\theta\in Irr(M)$ extends to $G$.
\end{prop}

We want to obtain a version of \cite[Lemma 2.4]{Brough3} for elements of prime power order, however it is not straight forward to construct an element $x$ of prime power order in ${\rm Sym}(n)$ such that 8 and every prime which divides $|{\rm Sym}(n)|$ also divides $|x^{\rm Sym}(n)|$.
Fortunately, for sporadic simple groups we do have the analogous result.

\begin{lm}\label{ListSporCases}
Let $S$ be a non-abelian sporadic simple group and assume there exists a prime $q$ such that $S$ does not have an irreducible character of $q$-defect zero.
\begin{enumerate}
\item There exists a prime power element $x$ whose conjugacy class $x^S$ is of size divisible by every prime dividing $S$ and by 4, and there exists $\theta\in Irr(S)$ which extends to $Aut(S)$ such that $\theta$ vanishes on $x^S$.
\item Let $p$ be a prime dividing the order of $S$.
Then there exists a $p'$-element $x$ of prime power order whose conjugacy class $x^S$ is of size divisible by $p$,  and there exists $\theta\in Irr(S)$ which extends to $Aut(S)$ such that $\theta$ vanishes on $x^S$.
\end{enumerate}

\begin{proof}
To prove $(1)$, the table below gives a pair $\{x_1,\theta_1\}$ satisfying the required conditions.

If a pair $\{x_1,\theta_1\}$ satisfies the conditions required for $(1)$, then it also satisfies the conditions required for $(2)$, unless $x_1$ turns out to have order divisible by $p$.
Thus to establish $(2)$ from $(1)$, it is enough to provide an additional pair $\{x_2,\theta_2\}$, such that if $x_1$ has order divisible by $p$ then $x_2$ has order not divisible by $p$. 
We cannot take the exact same list as in either \cite{Brough3}[Lemma 2.4] or \cite[Lemma 2.2]{DPSVan}, as the given elements were not all prime power elements.

The table below provides pairs $\{x_1,\theta_1\}$ and $\{x_2,\theta_2\}$ taken from \cite{Atlas}, for the required sporadic groups. 

\begin{center}
\begin{tabular}{lllll}
\hline
Group & Character $\theta_1$ & Class $x_1$ & Character $\theta_2$ & Class $x_2$\\
\hline
$M_{12}$ & $\chi_7$ & $3B$ & $\chi_7$ & $8A$\\
$M_{22}$ & $\chi_7$ & $8A$ & $\chi_2$ & $7A$\\
$M_{24}$ & $\chi_7$ & $4C$ & $\chi_5$ & $7A$\\
$J_2$ & $\chi_6$ & $3B$ & $\chi_{10}$ & $4B$\\
$HS$ & $\chi_7$ & $5C$ & $\chi_{16}$ & $4C$\\
$Suz$ & $\chi_3$ & $8B$ & $\chi_9$ & $3C$\\
$Ru$ & $\chi_{11}$ & $4D$ & $\chi_9$ & $5B$\\
$Co_1$ & $\chi_2$ & $4F$ & $\chi_2$ & $9B$\\
$Co_3$ & $\chi_6$ & $4B$ & $\chi_{10}$ & $5B$\\
$BM$ & $\chi_{20}$ & $4J$ & $\chi_{27}$ & $9B$\\
\end{tabular}
\end{center}

\end{proof}
\end{lm}

It remains to study the alternating groups.
We study ${\rm Alt}(n)$ for all $n\geq 7$, although in fact \cite[Corollary 2]{GODefectZero} yields some additional restrictions on $n$.
For $n\geq 7$ recall that Aut(Alt$(n)$)$\cong$Sym$(n)$.
As we are considering elements of prime power order it is enough to show the existence of such an element for a prime $l$ not equal to $2$ or $3$.
Then the simple group has a character of $l$-defect zero.

\begin{lm}\label{ListAltCases}
Let $S$ be a non-abelian simple group isomorphic to ${\rm Alt}(n)$ for $n\geq 7$ and assume there exists a prime $q$ such that $S$ does not have an irreducible character of $q$-defect zero.
\begin{enumerate}
\item There exists an $l$-element $x$ whose conjugacy class $x^S$ is of size divisible by 4 for some prime $l\ne 2$ or $3$.
\item Let $p$ be a prime dividing the order of $S$.
Then there exists an $l$-element $x$ whose conjugacy class $x^S$ is of size divisible by $p$,  for some prime $l\ne 2,3$ or $p$.
\end{enumerate}
\begin{proof}
To prove this statement, we first produce an $l$-element $x$ such that its conjugacy class size is divisible by 4 and every prime dividing ${\rm Alt}(n)$ except for $l$. 
In order to then obtain the second statement we can assume that the given prime $p$ is equal to $l$ for the example given to the first statement.
In this case it is then enough to produce another $l'$-element of prime power order with conjugacy class size divisible by $l$.

Let $l$ be the largest prime less than $n$ (i.e. $l$ is the largest prime dividing the order of ${\rm Alt}(n)$). 
Then as $n\geq 7$ it is clear that $l\geq 5$.
If $x$ is an $l$-cycle in ${\rm Alt}(n)$, then the size of its conjugacy class in ${\rm Sym}(n)$ is given by 
\[
\frac{n!}{l\cdot (n-l)!}=\frac{n(n-1)\dots (n-l+1)}{l}.
\]
As $l$ was chosen to be the largest prime less than $n$ it follows that both $4$ and every other prime divisor of ${\rm Alt}(n)$ not equal to $l$ divides the conjugacy class size of $x$ in ${\rm Alt}(n)$.
This completes the proof of the first claim.

Consider the second claim.
If the largest prime $l\leq n$ is not equal to $p$ then we are done.
Thus assume that $l=p$.
By using the verified Bertrand's-postulate \cite[Page 67]{NumThNag}, as $n\geq 7$, it follows that $l\leq n\leq 2l-1$.
Let $q$ be the second largest prime less than $n$, so $3<q\leq l\leq n$.
Let $k$ be a natural number such that $0\geq n-kq< l$.
If $k\geq l$, then $n-kq\leq n-pq\leq n-2p<0$ which is a contradiction.
Thus $k<l$.
Now let $x$ be a product of $k$ $q$-cycles.
It follows that the conjugacy class size of $x$ in ${\rm Sym}(n)$ is 
\[
\frac{n!}{q^k\cdot (k)!(n-kq)!}.
\]
However no term in the denominator of this fraction is divisible by $l$ and so $x$ is a $q$-element (i.e. a prime power $p'$-element) such that $p$ divides its conjugacy class size in ${\rm Alt}(n)$.
\end{proof}
\end{lm}

\section{The proofs}

\begin{thm}[Theorem 1]
 Let $G$ be a finite group and $p$ a prime divisor of $G$ such that if $q$ is any prime divisor of $G$, then $q$ does not divide $p-1$.
 Suppose that no conjugacy class size of a vanishing element of prime power order in $G$ is divisible by $p^2$. 
 Then $G$ is a soluble group.

 \begin{proof}

Suppose $G$ is chosen of minimal order satisfying the hypothesis of the theorem, but is not soluble.
By the same arguments as in \cite[Theorem A]{Brough3} it can be assumed $p=2$ and if $G$ has a proper normal subgroup $N$, then $G/N$ is soluble. 
Moreover a minimal normal subgroup $M\cong S_1 \times \cdots \times S_n$ is non-abelian.
If $S_i$ has a character of $q$-defect zero for all $q$ then as $S_i$ is non-soluble, \cite[Proposition]{LiTwoClaLen} implies that there is an element of prime power order with conjugacy class size divisble by $4$.
Meanwhile if $S_i$ is isomorphic to ${\rm Alt}(n)$ with $n\geq 7$ then by Lemma~\ref{ListAltCases}, there is an $l$-element for $l>3$ such that the conjugacy class size is disibe by $4$.
In both cases applying \cite[Lemma 2.2]{Brough3} shows $G$ has a prime power vanishing element with conjugacy class size divisible by $4$.
Hence it can be assumed $S_i$ is isomorphic to one of the sporadic groups given in Corollary~\ref{ListSimpleCases1}. 
In this case the same argument as in \cite[Theorem A]{Brough3} now using Lemma~\ref{ListSporCases} produces a vanishing element of $G$ with prime power order and 4 dividing its conjugacy class size.
\end{proof}
\end{thm}

\begin{thm}[Theorem 2]
 Let $G$ be a finite group and suppose that no conjugacy class size of a vanishing element of prime power order in $G$ is square free.
 Then $G$ is a super soluble group.
\begin{proof}
By combining Theorem 1 with the proof of \cite[Theorem B]{Brough3} this result now follows as the only elements considered are of prime power order.
\end{proof}
\end{thm}

\begin{thm}[Theorem 3]
 Let $G$ be a finite group and suppose a prime $p$ does not divide the size of any vanishing conjugacy class size $|x^G|$ for $x$ a $p'$-element of prime power order in $G$.
 Then $G$ has a normal $p$-complement.
 \begin{proof}
 Suppose $G$ is chosen of minimal order satisfying the hypothesis of the theorem, but does not have a normal $p$-complement.
By the same arguments used in the proof of \cite[Theorem C]{Brough3} we can conclude that $O_{p'}(G)=1$.
Let $M=S_1\times\dots\times S_k$ be a minimal normal subgroup of $G$, with each $S_i\cong S$ a simple group, then $p$ divides the order of $S$.
If $S$ is abelian the proof of \cite[Theorem C]{Brough3} shows that any vanishing $p'$-element of prime power order lies in $O_p(G)$.
Hence $G$ has a normal $p$-complement by \cite[Corollary C]{DPSVanOrd}.

Hence assume $S$ is non-abelian. 
First assume $S$ is not sporadic.
If $S$ has an irreducible character of $q$-defect zero for each prime $q$, then as $S\not\cong O_p(S)\times O_{p'}(S)$ it follows by \cite[Theorem 5]{LiWaWeNotesLenConjCl} that $S$ has a $p'$-element of prime power order such that $p$ divides its conjugacy class size. 
Moreover Lemma~\ref{ListAltCases} implies for $S$ of alternating type (on at least $7$ points) there exists a $p'$-element which has order a power of a prime $l>3$ and conjugacy class size divisible by $p$.
Thus \cite[Lemma 2.2]{Brough3} shows $G$ has a $p'$-element of prime power order which is vanishing and conjugacy class size divisible by $p$.
Finally assume $S$ is sporadic.
In this case the same argument as in \cite[Theorem C]{Brough3} now using Lemma~\ref{ListSporCases} produces a vanishing $p'$-element of $G$ with prime power order and $p$ dividing its conjugacy class size.
\end{proof}
\end{thm}

\section*{Acknowledgments}
The first author gratefully acknowledges financial support by the ERC Advanced Grant $291512$.
The research of the second author is supported by the National Natural Science Foundation of China (11301378).

\bibliographystyle{plain}
\bibliography{bibfile}

\end{document}